\documentclass{article}
\usepackage[utf8]{inputenc}
\usepackage{amsmath}
\usepackage{xcolor}
\usepackage{hyperref}

\newtheorem{note}{Note}
\newtheorem{alg}{Algorithm}
 \newtheorem{definition}{Definition}
\newtheorem{proposition}{Proposition}

\usepackage[colorinlistoftodos]{todonotes}

\begin{document}
\title{Experimental Mathematics Approach to Gauss Diagrams Realizability 
}

\author{
  A.~Khan\\
  University of Essex \\
  \texttt{ak20749@essex.ac.uk}
  \and
  A.~Lisitsa\\
  University of Liverpool \\
  \texttt{a.lisitsa@liverpool.ac.uk}
  \and
  A.~Vernitski\\
  University of Essex\\
  \texttt{asvern@essex.ac.uk}
  }

\maketitle

\begin{abstract}
A Gauss diagram (or, more generally, a chord diagram) consists of a circle and some chords 
inside it. 
Gauss diagrams are a well-established tool in the study of topology of knots and of planar and spherical curves. Not every Gauss diagram corresponds to a knot (or an immersed curve); if it does, it is called realizable. A classical question of computational topology asked by Gauss himself is which chords
diagrams are realizable. An answer was first discovered  in  the 1930s by Dehn, and since then many efficient algorithms for checking realizability of Gauss diagrams
have been developed. Recent studies 
in 
\cite{grinblat2018realizability,doi:10.1142/S0218216520500315} and \cite{doi:10.1142/S0218216519500159} formulated especially simple conditions related to realizability which are expressible in terms of parity of chords intersections. The simple form of these conditions opens an opportunity for experimental investigation of Gauss diagrams  using constraint satisfaction and related techniques. In this 
paper  we report on our experiments with Gauss diagrams of small sizes (up to 11 chords) using implementations of these conditions and other algorithms 
%in constraint satisfaction language MiniZinc and 
in logic programming language Prolog. In particular, we found a series of counterexamples showing that that realizability criteria 
%from [Biryukov, JKTR, Vol 28, No 1, 2019] 
in \cite{grinblat2018realizability,doi:10.1142/S0218216520500315} and \cite{doi:10.1142/S0218216519500159} are  not completely correct. 
\end{abstract}

\section{Introduction}	%) A SECTION HEADING

\begin{figure}[h!]
\centering
\includegraphics[scale=0.4]{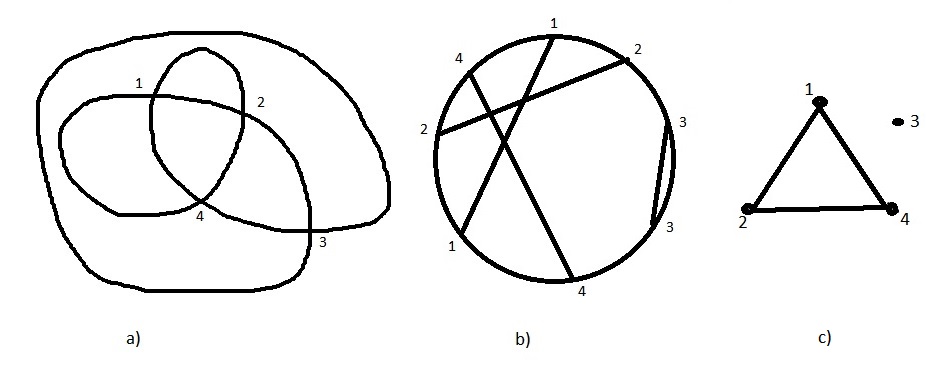}
\caption{Example of a) a planar curve; b) its Gauss diagram and c) its interlacement graph. The corresponding Gauss word is {\bf 12334124} }
\label{fig:gauss-ex}
\end{figure}

Gauss words and Gauss diagrams are natural combinatorial encodings of generic smooth closed planar curves.
For a closed planar curve its Gauss code (word) can be obtained by first assigning different symbols to all intersection points and then 
by writing down the symbols of intersection points passed over during traversal of the curve, 
starting and finishing at the chosen non-intersection  point, see (Fig~\ref{fig:gauss-ex}, a). 
The Gauss code of a closed curve with $n$ intersection points 
has a length $2n$ and it is a \emph{double occurrence} word, that is in which each symbol occurs twice. 
With  any \emph{double occurrence} word 
one can further associate its  chord (or Gauss) diagram which consists of a circle with all symbols of the word placed at the points of the circle following its clockwise traversal. The chords of the diagram link the points labelled by the same symbol, (Fig.\ref{fig:gauss-ex},b). 
Not every double occurrence word and  its Gauss diagram correspond to  (can be obtained from)  a closed planar curve; if they do, they are called \emph{realizable}.  

\begin{figure}[h!]
\centering
\includegraphics[scale=0.5]{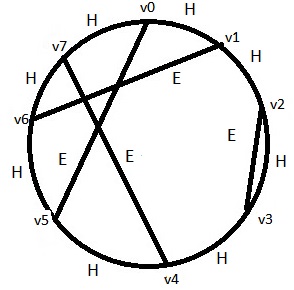}
\caption{Diagram graph of the diagram from Fig 1, b)}
\label{fig:gauss-graph}
\end{figure}

A classical question of computational topology asked by Gauss \cite{Gauss} is which of the chords diagrams are realizable. 
Gauss himself  gave  a partial answer in terms of a necessary condition - \emph{in a realizable diagram  any chord intersects even number  
of other chords} 
Additional necessary condition was provided by Nagy \cite{n-utpg-27} and a full answer to Gauss' question was discovered in the 1930s by Dehn \cite{dehn1936}. Since then many efficient criteria and  algorithms for checking the realizability of Gauss diagrams 
have been developed
\cite{10.2307/2037443,FRANCIS1969331, DEFRAYSSEIX198129,lovasz1976,rosenstiehl:hal-00259712, rosenstiehl:hal-00259721,DBLP:journals/jal/RosentiehlT84,10.1007/3-540-63938-1_65,DOWKER198319,vena2018topological,ce-pp2-96,cw-prpmg-94}. 
As computational complexity concerned recognition of realizability of Gauss diagrams or words can be done in \emph{polynomial time} 
and even in linear time \cite{DBLP:journals/jal/RosentiehlT84}. Furthermore as it is shown in \cite{10.1007/978-3-642-21254-3_29}  it can be done with \emph{logarithmic space} complexity. All  realizability criteria and detection algorithms in the papers cited above  can be expressed in terms of a \emph{diagram graph} $(V,E,H)$ with two types of edges representing a chord diagram. Vertices of the graph are points of the diagrams, one type $E$ of edges corresponds to chords, the edges of another type  $H$
form a  
cycle covering all vertices in $V$, see Fig.~\ref{fig:gauss-graph}.

With every Gauss diagram one may further associate an \emph{interlacement} graph whose vertices are the chords of the diagram and two vertices are adjacent if an only if the corresponding chords intersect, see (Fig.~\ref{fig:gauss-ex}, c).
The interlacement graph of a diagram embodies an essential information about its realiziability. Furthermore, as it was shown in  
\cite{rosenstiehl:hal-00259712,10.1007/3-540-63938-1_65,DBLP:journals/dm/ShtyllaTZ09}  
the realizability  of Gauss diagrams  
can be decided based on 
their 
interlacement graphs only.  The criteria are more abstract then those defined in terms of diagram graphs, as different (non-isomorphic) diagram graphs may have isomorphic corresponding interlacement graphs. At the same time criteria from  \cite{rosenstiehl:hal-00259712,10.1007/3-540-63938-1_65,DBLP:journals/dm/ShtyllaTZ09}  appear to be computationally more costly. Some account of this can be given in terms of definability and descriptiive complexity~\cite{DBLP:books/daglib/0095988}.    The formal definition of these conditions  %\cite{DBLP:journals/dm/ShtyllaTZ09} 
can be presented in  first-order logic (FO) extended by  a parity quantifier (${\rm MOD_{2}}$)~\cite{NURMONEN200062} and monadic second-order quantification (MSO). 
%\cite{definability}. 
While parity quantifier evaluation can be implemented very efficiently, second-order quantification goes over all subsets of vertices.      Hence checking these conditions by a straightforward algorithm requires \emph{exponential time}.     
More recent studies \cite{grinblat2018realizability,doi:10.1142/S0218216520500315,doi:10.1142/S0218216519500159} formulated  
even simpler realizability conditions 
expressible in terms of the interlacement graphs.  
In particular \cite{grinblat2018realizability,doi:10.1142/S0218216520500315} proposed a criteria 
based on 
parity conditions and \emph{reduction} operation 
on 
the graph, while   \cite{doi:10.1142/S0218216519500159} reformulated it in a 
simple form, expressed  in  ${\rm FO +MOD_{2}}$.
Further simple analysis shows that the original conditions from \cite{grinblat2018realizability} can be expressed in ${\rm FO +MOD_{2}}$  too. Hence both variants of criteria  
can be tested 
in polynomial time and logarithmic space.
Furthermore the simple and declarative  form of these conditions opens an opportunity for experimental investigation of Gauss diagrams  using constraint satisfaction \cite{books/daglib/0076790,10.5555/1207782} and related computational techniques.   
We took this opportunity as a starting point of our research in experimental mathematics with the primary aim to enumerate various classes of Gauss diagrams. 
We have implemented the declarative realizability conditions and algorithms from \cite{doi:10.1142/S0218216519500159}, \cite{doi:10.1142/S0218216520500315} and    
\cite{DBLP:journals/dm/ShtyllaTZ09} alongside with the classical algorithm from \cite{dehn1936} (for the latter see detailed exposition in e.g. \cite{KAUFFMAN1999663}).  We have enumerated several interesting classes of Gauss diagrams up to sizes 10-12 (number of chords), confirming and expanding on known enumerations (for the latter see e.g. \cite{Coquereaux:2015eia} and references therein). The detailed presentation of this will be given in 
an extended 
version of 
this 
paper. In the present note  we confine ourselves with the presentation of the most surprising discovery, that is the realiziability conditions in \cite{doi:10.1142/S0218216519500159}, \cite{doi:10.1142/S0218216520500315} are 
not completely correct. 
We present a minimal counterexample to both conditions of size 9 and enumerate all countexamples of sizes 10 and 11. We further confirm experimentally up to the size 11 that criteria from \cite{doi:10.1142/S0218216519500159} and  \cite{doi:10.1142/S0218216520500315} are equivalent,  and  the condition from \cite{DBLP:journals/dm/ShtyllaTZ09} is correct, i.e equivalent to that from \cite{dehn1936}.

\section{Preliminaries}

In this section we formally discuss the definitions used in our paper.

\begin{definition}\cite{c-cic-91}  
Immersed into a surface G curve is a mapping $\gamma: S^{1} \rightarrow G$.
Two curves $\gamma_{1}$ and $\gamma_{2}$ 
are geotopic~\cite{c-cic-91}, or just equivalent  if there is a homeomorphism 
$h: G \rightarrow G$ such that $h\gamma_{1} = \gamma_{2}$.   
\end{definition} 

In this paper we are only concerned with the curves immersed into plane $R^{2}$ or sphere $S^{2}$.

\begin{definition}~\label{g(w)}
A word is called double occurrence 
if each letter appears exactly twice in it. 
With every double occurrence word $w = a_{0}a_{1}\ldots a_{2n-1}$, one can associate a Gauss diagram $g(w) = \langle V,E,H\rangle$ , that is an undirected  graph with two types of edges. Here the set of vertices $V =\{ i \mid 0 \le i \le 2n-1 \}$, 
the set of edges called chords $E= \{(i,j)\mid i < j \; \& \; a_{i} = a_{j}\}$, and set of edges $H = \{(i,j)\mid j = i + 1\; mod\; 2n  \}$ which form a cycle .
\end{definition}

\begin{definition} 
An abstract Gauss diagram is a finite undirected graph $\langle V,E,H\rangle$ with the set of vertices $V$ and two types of edges $E$ and $H$   such that $E$ is a perfect matching on $V$ (i.e. each vertex has exactly one neighbour in $E$), and $H$ forms a simple  cycle covering all vertices in $V$. We will occasionally omit the qualifier word  ``abstract'' for Gauss diagrams when it is clear from the context and/or does not lead to the confusion.     
\end{definition}

\begin{definition}
Given a Gauss diagram $g = \langle V,E,H \rangle$, the interlacement graph $I_{g}$ of $g$ is an undirected  graph $\langle V_{I},E_{I}\rangle$, where the set of vertices $V_{I}$ is the set of edges $E$ of $g$, and 
$((i,j),(k,l)) \in E_{I}$ iff the edges $(i,j)$ and $(k,l)$ are interlaced in $g$, that is either $i < k < j < l$, or $k < i < l < j$ holds true.      \end{definition}

\begin{definition} 
A Gauss diagram $g$ is called \emph{prime}
if and only if its interlacement graph $I_{g}$  is connected. \end{definition}

The  equivalence of immersed curves with prime Gauss diagrams  is captured by \emph{isomorphism}  of Gauss words and Gauss diagrams. 
Two double occurrences words over the same alphabet are  isomorphic  if one can be obtained from the other by  any combination of cyclic shifts and reversing the word \cite{c-cic-91}.  
For immersed curves $\gamma_{1}$ and $\gamma_{2}$, such that their Gauss diagrams are prime,  both  corresponding  Gauss words and Gauss diagrams are isomorphic if and only if $\gamma_{1}$ and $\gamma_{2}$ are equivalent,  see e.g.  \cite{chmutov-2006}.

A double occurrence word $w$ is called \emph{realizable}  
iff there exists an immersed in a plane (or sphere $S^{2}$\footnote{when presenting criteria for realizability of Gauss words and diagrams, various authors refer to immersions either to plane, or sphere $S^{2}$. It is easy to see, by considering the  stereographic projection of $S^{2}$ to $R^{2}$, that it does not make a difference })  curve $\gamma$, such that 
$w$ is associated with $\gamma$. A Gauss diagram $g$ is \emph{realizable} iff there is a realizable  word $w$ such that $g$ is  isomorphic to $g(w)$. 

Not all double occurrences words and Gauss diagram are realizable. In this paper we mainly focus on investigation of the realizability criteria of prime Gauss diagrams expressed solely in terms of their interlacement graphs. 

One necessary realizability condition  proposed by Gauss himself can be formulated in terms of interlacement graphs  as follows: 

\begin{quote}
{\bf C1:} \emph{If a Gauss diagram $g$ is realizable then in its interlacement graphs $I_{g}$ every node has an even degree. }
\end{quote}

\section{Realizability criteria based on interlacement graphs} 

In this section we overview known realiziability criteria for Gauss diagrams, which are defined in terms of the interlacement graphs.

\subsection{R-conditions} 

For the first time
relizability criteria for Gauss words and diagrams expressed solely in terms of interlacement graphs  have 
been established by Rozentiehl in 
\cite{rosenstiehl:hal-00259712}. Following  \cite{10.1007/3-540-63938-1_65} Rozentiehl's conditions can be formulated as follows. 
A Gauss diagram $g$ is realizable iff its interlacement graph $I_{g}$ has the following properties: 

\begin{itemize}
    \item $I_{g} = (V,E)$ is eulerian, that is each vertex has an even degree; \hspace*{8mm} {\bf (C1)}
    \item there is a subset of vertices $A \subset V$ of $I_{g}$ such that the following two conditions are equivalent for any two vertices $u$ and $v$:  
    \begin{itemize}
    \item the vertices u and v have an odd number of common neighbours, 
    \item the vertices u and v are neighbours and either both are in $A$ or neither is in $A$  
    \end{itemize}
\end{itemize}

\subsection{STZ-conditions} 

Shtylla, Traldi and Zulli 
in \cite{DBLP:journals/dm/ShtyllaTZ09}, while keeping {\bf C1} condition,  
presented an algebraic re-formulation of the second Rosentiehl's condition: 

\begin{itemize}
\item For its adjacency matrix $M_{g}$ there is a 
    diagonal matrix $\Lambda_{g}$ such that $M_{g}+\Lambda_{g}$ is an idempotent matrix (all matrices are considered over $GF(2)$, the finite field of two elements). %\todo{not defined whats Gf(2)}.
\end{itemize}

\subsection{GL-conditions}
Further criteria for realizability expressible in terms of interlacement graph  have been presented by Grinblat and Lopatkin in  \cite{grinblat2018realizability,doi:10.1142/S0218216520500315}. A prime Gauss diagram $g$ is realizable if and only if the following conditions hold true: 

\begin{itemize}
    \item In its interlacement graph $I_{g}$ each pair of non-neighbouring vertices has an \emph{even} number of common neighbours (possibly, zero).  $\hspace*{2cm}$   ({\bf C2})
    \item The above condition {\bf C2} holds true for the reduced graph $I_{g}/v$ for each vertex $v$ of $I_{g}$.   
\end{itemize}

For a vertex $v$ of $I_{g} = \langle V,E\rangle$ the reduced graph $I_{g}/v = \langle V',E' \rangle$ is defined as follows.
$V' = V -\{v\}$ and  $E' = \{ (v_{1},v_{2}) \in E |  ((v,v_{1}) \not\in E)\lor (v,v_{2}) \not\in E))  \} \cup  \{ (v_{1},v_{2}) \mid (v,v_{1}) \in E \;\& \; (v,v_{2}) \in E \;\&\; (v_{1},v_{2}) \not\in E\}$. See also further discussion of this operation in terms of intersecting chords in Section 2 of   
\cite{doi:10.1142/S0218216519500159}.

\subsection{B-conditions} 

Biryukov in \cite{doi:10.1142/S0218216519500159} presents even simpler 
than GL-conditions realizability criteria  explicitly expressed in terms of interlacement graph: 
A prime Gauss diagram $g$ is realizable  if and only if 
the following conditions for 
$I_{g} = \langle V,E\rangle$ 
hold true: 

\begin{itemize}
    \item It  satisfies both above {\bf C1} and {\bf C2} conditions, which amounts to  \emph{strong parity condition} from  \cite{doi:10.1142/S0218216519500159};
    \item  For any three pairwise connected vertices  $a, b, c \in V$   the sum of the number of vertices  adjacent to  $a$, but not adjacent to  $b$ nor $c$, and the number of vertices adjacent to  b and c, but not adjacent to  $a$, is even.   \hspace*{2cm}       {\bf(B3)}
\end{itemize}

\subsection{Remarks on definability and complexity} 
$R$- and $STZ$- conditions use explicit second order-quantifier ``there exists a subset A'' and can be defined over $I_{g}$ in monadic second order logic extended by a parity quantifier (${\rm MSO+MOD_{2}}$). That  implies the exponential time upper bound for computational complexity of their checking. B-conditions on the other hand don't use second-order quantification and can be  defined in first-order logic extended by parity quantifier (${\rm FO+MOD_{2}}$). 
For example,  the definition for  {\bf(B3)} is: 
\[ \forall x, y, z ((E(x,y) \land E(x,z) \land E(y,z) \Rightarrow \] 
\[ \exists^{even} w ((E(x,w) \land \neg E(y,w) \land \neg E(z,w) \land(y \not = w) \land (z \not = w)) \] \[ \lor (\neg E(x,w) \land E(y,w) \land E(z,w) \land (x \not =w)))) \]

Here $E(\_,\_)$ refers to (interpreted by) the edges of the interlacement graph and  
$\exists^{even} w$ is parity quantifier meaning 
``there exists an \emph{even} number of elements/valuations of $w$  such that $\ldots$''. 
${\bf C2}$ condition  can be  defined in a similar way. Furthermore,  using \emph{locality} property of  the reduction operation, the second GL condition can also be defined in 
${\rm FO+MOD_{2}}$.
Thus, both B and GL conditions can be checked with logarithmic space (and polynomial time)  complexity \cite{DBLP:books/daglib/0095988}.

\section{Lintels}

We aim to implement efficient algorithms for generating Gauss diagrams, checking their (non-)equivalence, listing all non-equivalent diagrams, satisfying combinations of criteria and  checking the realiziability criteria. As a starting point we take Definition~\ref{g(w)}, and consider the generation of Gauss diagrams with a sets of vertices $V =\{ i \mid 0 \le i \le 2n-1 \}$ for some $n \ge 2$ satisfying an additional  
constraint, corresponding to the condition C1. For that we introduce a concept of \emph{lintel},
and useful algorithms working on lintels.

Let $n$ be a positive integer. A \emph{lintel} is an $n$-tuple of pairs of numbers $((a^1_1, a^2_1), \dots,$ $(a^1_n, a^2_n))$ such that 1) the set of numbers in the lintel $\{a^1_i, a^2_i, \dots, a^1_n, a^2_n\}$ is equal to the set $\{ 0, 1, \dots, 2n-1 \}$, and 2) in each pair $(a^1_i, a^2_i)$ the difference $a^1_i - a^2_i$ is odd. The number $n$ will be called the size of the lintel. Each pair $(a^1_i, a^2_i)$ will be called a chord of the lintel.

For example, $L = ((4, 7), (8, 1), (5, 0), (2, 9), (6, 3))$ is a lintel of size $5$. For better readability, below we will write the two numbers in each chord under one another; therefore, the lintel $L$ can be written as $\begin{pmatrix}
4 & 8 & 5 & 2 & 6\\
7 & 1 & 0 & 9 & 3
\end{pmatrix}$.

Let us say that two lintels are \emph{strongly equivalent} if they can be transformed into one another by steps of the following two types: 
\begin{enumerate}
    \item Swapping the positions of two numbers in a chord; for example, the lintel $L$ in the example above is strongly equivalent to $\begin{pmatrix}
4 & 8 & 0 & 2 & 6\\
7 & 1 & 5 & 9 & 3
\end{pmatrix}$;
    \item Swapping the positions of chords in the list; for example,  the lintel $L$ in the example above is strongly equivalent to $\begin{pmatrix}
2 & 8 & 5 & 4 & 6\\
9 & 1 & 0 & 7 & 3
\end{pmatrix}$.
\end{enumerate}

Among the lintels that are strongly equivalent to each other, one can point at one which can serve as the canonical lintel of this class of lintels. To define it, let us introduce a kind of lexicographic order on lintels. Consider two lintels of size $n$, $L = \begin{pmatrix}
a^1_1 & \dots & a^1_n\\
a^2_1 & \dots & a^2_n
\end{pmatrix}$ and $M = \begin{pmatrix}
b^1_1 & \dots & b^1_n\\
b^2_1 & \dots & b^2_n
\end{pmatrix}$. We will say that $L$ is less than $M$ if there is a position $i$ such that $a^1_j = b^1_j$ and $a^2_j = b^2_j$ for all $j < i$, and either $a^1_i < b^1_i$, or $a^1_i = b^1_i$ and $a^2_i < b^2_i$. We will call this order the \emph{L-order}. The lintel which is the least one relative to the L-order among the lintels that are strongly equivalent to each other will be treated as the canonical representative of this class of lintels. There is also another way to define this lintel. Let us say that a lintel is a \emph{sorted lintel} if each pair in it is sorted, and first elements of pairs are sorted, that is, for each $i$ we have $a^1_i < a^2_i$, and for all $i < j$ we have $a^1_i < a^1_j$. In each class of lintels that are strongly equivalent to each other, there is exactly one sorted lintel, and it is the one that is the least one relative to the L-order. For a given lintel $L$, one can produce the sorted lintel corresponding to $L$ using the following simple algorithm.

\begin{alg}

1. Sort the two numbers in each chord.
\\2. Sort all the chords by the values of their first entry.
\end{alg}

Let us say that two lintels are \emph{equivalent} if they can be transformed into one another by steps of the following types: 
\begin{enumerate}
    \item the two types of steps defining strongly equivalent lintels
    \item shifting the value of each entry of the lintel cyclically modulo $2n$, that is, replacing a lintel $\begin{pmatrix}
a^1_1 & \dots & a^1_n\\
a^2_1 & \dots & a^2_n
\end{pmatrix}$ by a lintel $\begin{pmatrix}
a^1_1 + s & \dots & a^1_n + s\\
a^2_1 + s & \dots & a^2_n + s
\end{pmatrix}$ for some value $s$, with operations performed in the arithmetic modulo $2n$;
    \item inverting the value of each entry of the lintel modulo $2n$, that is, replacing a lintel $\begin{pmatrix}
a^1_1 & \dots & a^1_n\\
a^2_1 & \dots & a^2_n
\end{pmatrix}$ by the lintel $\begin{pmatrix}
- a^1_1 & \dots & - a^1_n\\
- a^2_1 & \dots & - a^2_n
\end{pmatrix}$, with operations performed in the arithmetic modulo $2n$;
\end{enumerate}

Now we are ready to define a canonical lintel. 
Namely, the lintel which is the least one relative to the L-order among the lintels that are equivalent to each other will be treated as the canonical representative of this class of lintels 
We shall call such lintels \emph{Lyndon lintels}. This construction is similar in certain aspects to Lyndon words in algebra, hence the name
Obviously, a Lyndon lintel is a sorted lintel, but not every sorted lintel is a Lyndon lintel. For a given lintel $L$, one can produce the sorted lintel corresponding to $L$ using the following algorithm.

\begin{alg}
1. Produce the sorted lintel $M$ corresponding to $L$.
\\2. Applying cyclic shifts to $M$ for all values $s = 1, \dots, 2n-1$, produce lintels $L_s$.
\\3. Produce the sorted lintels $M_s$ corresponding to $L_s$.
\\4. Applying the inverting step to $L$, produce the lintel $L'$.
\\5. Produce the sorted lintel $M'$ corresponding to $L'$.
\\6. Applying cyclic shifts to $M'$ for all values $s = 1, \dots, 2n-1$, produce lintels $L'_s$.
\\7. Produce the sorted lintels $M'_s$ corresponding to $L^-_s$.
\\8. Among the sorted lintels $M, M', L_s, L'_s$, choose the least one relative to the L-order.
\end{alg}

\begin{proposition}\label{lin-per} 
The number of equivalence classes for strong equivalence relation on lintels of size $n$ is $n!$
\end{proposition} 
{\bf Proof}
Let $I_{2n} = \{0, \ldots, 2n-1\}$. Denoted by $P_{2n}$ the set of all even-odd matchings of $I_{2n}$, that is partitions of $I_{2n}$ into subsets of size 2, each containing one even and one odd number.
The obvious consequence of definitions is that  elements of $P_{2n}$ are in bijective correspondence with the set of   all equivalence classes  for strong equivalence relation on lintels of size $n$. Now for symmetric group $S_{n}$ we define a mapping $ \beta: S_{n} \rightarrow P_{2n}$ as follows. 
For  $\sigma \in S_{n}$,  $\sigma:[1..n]\rightarrow[1..n]$  $\beta(\sigma) =  \{\{2*i-1,2*\sigma(i)-2\}| i \in [1..n]\}$. Straightforward check shows that $\beta$ is a bijection between $S_{n}$ and $P_{2n}$.

\begin{note}
The fact $|P_{2n}| = n!$ is mentioned (without the reference to the proof) in 
OEIS article A000142 \cite{oeisA000142}  with an attribution to  David Callan, Mar 30 2007. 
We use  explicit simple bijection between $P_{2n}$ and  $S_{n}$ 
for the efficient generation of 
lintels. 
\end{note}

\section{Implementation}
We have implemented the algorithms suite \verb"Gauss-Lint" for Gauss diagram generation, equivalence checking,  canonization, and realizability  conditions using the concept of lintels. The prototype implementation in logic programming language SWI-Prolog \cite{wielemaker:2011:tplp}  is available at \cite{abdullah_khan_2021_4574590}.
 
Here we give a short overview of the implementation while  more detailed presentation is deferred to the extended version of this paper. 

\subsection{Essential feature}

As we use lintel concept/data structure, all generated Gauss diagrams in our implementation satisfy condition ${\bf C1}$. 

\subsection{Data structures}
We utilise a  list of lists representation of lintels in Prolog. 
As an example, Prolog list \verb"[[0,3],[1,4],[2,5]]"
 represents a lintel $\begin{pmatrix}
0 & 1 & 2 \\
3 & 4 & 5 
\end{pmatrix}$.

\subsection{Dynamic facts} 
We use dynamic facts to store intermediate and final results of computations. 
E.g. a dynamic fact ``lintel(Size,L)'' stores a generated canonical lintel satisfying some conditions (depending on the user's choice). 

\subsection{Lintel generation and canonization}

The canonical lintel generation is based on efficient built-in Prolog predicate \verb"permutation(X,Y)" for the generation of all permutations of a given list (assigned to) $X$ and the use of bijection 
$\beta: S_{n} \rightarrow P_{2n}$ (Section 4). 

The generic generation process proceeds as follows. 
Once a permutation $\sigma$ is generated, the bijection $\beta$ is applied to it, and then the Algorithm 2 is 
applied to $\beta(\sigma)$. The resulted canonical lintel then checked on 
\begin{itemize}
\item whether it satisfies chosen combination of conditions;  
\item whether it has not  been stored yet as a dynamic fact.
\end{itemize}
If both conditions are satisfied then the lintel is added as a dynamic fact, and the process is backtracked  to the generation of the new permutation. When no new permutations are  left, the database of dynamic facts provides a list of all canonical non-equivalent lintels of given size satisfying chosen combination of properties.   

\subsection{Properties/conditions checking}
 We have implemented the checking of following properties/conditions of Gauss diagrams presented in lintel format: 
 primality, C2, STZ, B, GL. As mentioned before, C1 is satisfied automatically for diagrams presented by lintels. We have also implemented a classical  algorithm (CA) for realiziability checking from \cite{dehn1936,KAUFFMAN1999663}.      

\section{Experiments and  Counterexamples} 

Using our implementation we have conducted various experiments and enumerated the classes of non-equivalent Gauss diagrams, satisfying various    
combinations of properties.  

\begin{table}[h]
 \centering
    \label{tab:table1}
    \begin{tabular}{|l|c|c|c|c|c|c|c|c|c|c|} % <-- 
   \hline 
     & 3 & 4 & 5 & 6 & 7 & 8 & 9 & 10 & 11 
      %& 12 
     \\
      \hline
      Realisability (CA) & 1 & 1 & 2 & 3 & 10 & 27 & 101 & 364 & 1610 
      %& 7202
      \\
      STZ & 1 & 1 & 2 & 3 & 10 & 27 & 101 & 364 & 1610 
      %& $\ldots$ 
      \\
      B  & 1 & 1 & 2 & 3 & 10 & 27 & {\bf 102} & {\bf 370} & {\bf 1646} 
      %& {\bf 7437} 
      \\
      GL & 1 & 1 & 2 & 3 & 10 & 27 & {\bf 102} & {\bf 370} & {\bf 1646} 
      %& $\ldots$ 
      \\
      \hline
    \end{tabular}
   \caption{The number of non equivalent Gauss diagrams of sizes  = 3, \ldots, 11, satisfying various realizability conditions}
\end{table}

Table 1 shows the numbers of non-equivalent Gauss diagrams of sizes $3\ldots 11$ satisfying different realizability conditions, generated by our program.  The first two lines, found also in OEIS article A264759 \cite{oeisA264759}, together with the fact that generated corresponding lists of Gauss diagrams are the same,   verify that indeed, STZ are correct realizability conditions up to the size 11.  

The last two lines show an interesting phenomenon.  Up to the size 8 B and GL conditions are satisfied by the same diagrams, as  CA and STZ are.  For the size 9, however, there is one (= $102 - 101$)  up to equivalence Gauss diagram such that it satisfies both B and GL conditions, but still is not realizable. For sizes 10 and 11 the numbers of such diagrams are 6 and 36\footnote{ For size 12, the number of such diagrams is $235 = 7437-7202$. The analysis of size 12 required an additional, more efficient method of generating Gauss diagrams, which will be described in the extended version of this paper 
}, respectively. 

In summary, CA and STZ, respectively  B and GL, are pairwise equivalent conditions up to the size 11. B and GL, however are not equivalent to CA or STZ, starting from size 9. Hence B and GL are not correct realizability criteria despite the claims in \cite{doi:10.1142/S0218216519500159} and \cite{grinblat2018realizability,doi:10.1142/S0218216520500315}, respectively.     

\subsection*{Counterexample, Size=9} 

\begin{figure}[h!]
\centering
\includegraphics[scale=0.5]{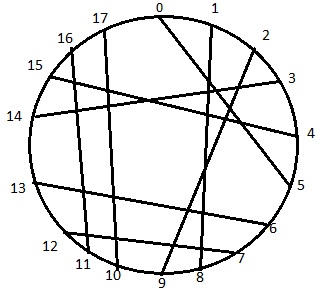}
\caption{Minimal counterexample, n=9}
\label{fig:min}
\end{figure}

The lintel 
\begin{verbatim}
[[0,5],[1,8],[2,9],[3,14],[4,15],[6,13],[7,12],[10,17],[11,16]]
\end{verbatim}
represents  a minimal and single up to equivalence non-realizable Gauss diagram of size 9, satisfying both B and GL conditions, Fig~\ref{fig:min}.  

\subsection*{Counterexamples, Size=10}
\begin{figure}[h!]
\centering
\includegraphics[scale=0.4]{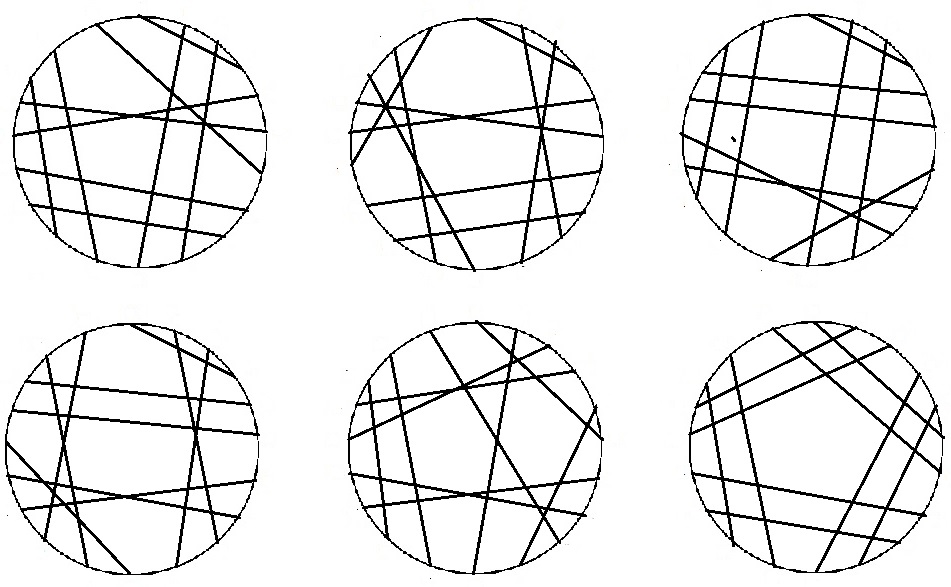}
\caption{All counterexamples for n=10}
\label{fig:counter10}
\end{figure}

The following lintels represent all six up to equivalence non-realizable Gauss diagrams of size 10 satisfying both B and GL conditions, Fig~\ref{fig:counter10}. 

\begin{verbatim}
[[0,3],[1,10],[2,9],[4,15],[5,16],[6,19],[7,14],[8,13],[11,18],[12,17]]
[[0,3],[1,8],[2,9],[4,15],[5,16],[6,13],[7,12],[10,17],[11,18],[14,19]]
[[0,3],[1,10],[2,9],[4,17],[5,16],[6,11],[7,14],[8,15],[12,19],[13,18]].
[[0,3],[1,8],[2,9],[4,17],[5,16],[6,13],[7,14],[10,15],[11,18],[12,19]].
[[0,5],[1,10],[2,15],[3,16],[4,9],[6,13],[7,14],[8,19],[11,18],[12,17]].
[[0,5],[1,16],[2,15],[3,10],[4,9],[6,19],[7,14],[8,13],[11,18],[12,17]].
\end{verbatim}

\section{Conclusion}

We have presented an approach to experimental investigation of Gauss realizability conditions based on the concept of lintel. We have implemented the algorithms for Gauss diagrams generation, canonization and checking recently proposed declarative criteria for Gauss realizability. We reported on the experiments with our implementations and presented  counterexamples for the recently published 
simple realizability criteria. Based on that we conjecture that 
realizability of Gauss diagrams can not be defined in first-order logic with parity quantifier over their interlacement graphs.

\section*{Acknowledgments}
This work was supported by the Leverhulme Trust Research Project Grant RPG-2019-313.

\bibliographystyle{alpha}
\bibliography{sample}

\end{document}